\documentstyle[12pt]{article}
\begin{document}
\title{\bf The Baer Radical Of Generalized Matrix Rings }
\author{\em Shouchuan Zhang \\
Department  of Mathematics, Hunan
University,\\
Changsha  410082, P.R.China, E-mail:Zhangsc@hunu.edu.cn}
\date{}
\newtheorem{Theorem}{Theorem}[section]
\newtheorem{Proposition}{Proposition}[section]
\newtheorem{Definition}{Definition}[section]
\newtheorem{Corollary}{Corollary}[section]
\newtheorem{Lemma}{Lemma}[section]
\newtheorem{Example}{Example}[section]
\maketitle \addtocounter{section}{-1}
\begin {abstract}
In this paper, we introduce a new concept of generalized matrix
rings and build up the general theory of radicals for g.m.rings.
Meantime,  we obtain
$$\bar{r}_b(A)=g.m.r_b(A)=\sum\{r_b(A_{ij})\mid
i, j\in I\}=r_b(A)$$ Key words:$\Gamma$-ring,  small additive
category,  generalized matrix ring,  radical,  directed graph.¦Ì\\
AMS(1989) Subject classifications.16A21 16A64, 16A66, 16A78
\end {abstract}

\section {Introduction}
Let $I$ be a set.  For any  $ i, j, k\in I$, $(A_{ij}, +)$ is an
additive abelian group and there exists a map $\mu_{ijk}$ from the
product set $A_{ij}\times A_{jk}$ of $A_{ij}$ and $A_{jk}$into
$A_{i k}$(written $\mu_{ijk}(x, y)=xy$) such that the following
conditions are satisfied.

$(i) (x+y)z=xz+yz, w(x+y)=wx+wy.$

$(ii)w(xz)=(wx)z.$\\
for any $x, y\in A_{ij}, z \in A_{jk}, w\in A_{li}$.  We call
$\{A_{ij}\mid i, j\in I\} $ a $\Gamma_I$-system.

An additive category ${\cal C}$ is a small additive category if
ob${\cal C}$ is a set.

If ${\cal C}$ is a small additive  category and $I=ob{\cal C},
A_{ij}=Hom_{\cal C} (j,  i)$ for any $i, j\in I$, then we easily
show that $\{ A_{ij}\mid i, j\in  I\}$ is a $\Gamma_I$-system.

$D$ is a directed graph (simple graph).   Let $I$ denote the
vertex set of $D$ and $x=(x_1, x_2, \cdots ,  x_n)$ a path from
$x_1$ to $x_n$ via vertexes  $x_2, x_3, \cdots ,  x_{(n-1)}$, If $
x=(x_1, x_2, \cdots , x_n)$ and $y=(y_1, y_2, \cdots ,  y_m)$ are
paths of $D$ with $x_n=y_1$, we define the multiplication of $x$
and $y$
$$xy=( x_1, x_2, \cdots , x_n,   y_2, \cdots ,  y_m)$$

Let $A_{ij}$ denote a free abelian group generated by all paths
from $i$ to $j$, where $i, j\in I$.  We can define a map from $
A_{ij}\times A_{jk}$ to $A_{ik}$ as follows.
$$(\sum_s(n_sa^{(s)}))(\sum_sm_tb^{(t)})=
\sum_s\sum_tn_sm_ta^{(s)}b^{(t)}$$where $a^{(s)}\in A_{ij},
b^{(t)}\in A_{jk}, s, t\in {\bf Z}^+, n_s, m_t\in {\bf Z}$.  We
easily show that $\{A_{ij}\mid i, j\in I\}$ is a $\Gamma
_I$-system.

Let $\{A_{ij}\mid i, j\in I\}$ be a $\Gamma _I$-system and $A$ the
external direct sum (if $x\in A$, then $x_{ij}=0$ for all but a
finite number of $i, j\in I$).  We define the multiplication in $A$
as follows.  $$xy=\{\sum_kx_{ik}y_{kj}\}$$ for any $x=\{x_{ij}\},
y=\{y_{ij}\}\in A$.

It is easy to show $A$ is a ring.   We call $A$ a generalized
matrix ring,  written $A=\sum.  \{A_{ij}\mid i, j\in I\}$.  Every
element in $A$ is called a generalized matrix.

Therefore,  every generalized matrix ring is a ring.   Conversely,
every ring $A$ can also be considered as a generalized matrix ring
$A$ with
$$A=\sum \{A_{ij}\mid i, j=1\}$$

In this paper,  $A$ denotes a generalized matrix ring.  $xE(i, j)$
denotes a generalized matrix having a lone $x$ as its $(i,
j)$-entry all other entries $0$.  ${\bf Z}^+$ denotes the natural
number set.   ${\bf Z }$ denotes the integer number set.

Every generalized matrix ring can be written g.m.ring in short.

Shaoxue Liu obtained that the relations between the radical of
additive category $A$ and the radical of ring $A_{ii}$ in [3].

In this paper we obtain that the relations between the Baer
radical of g.m.ring $A$ and the Baer radical of $A_{ji}$-ring
$A_{ij}$.
\section {The Basic Concept}
In this section we omit all of the proofs because they are
trivial.
\begin {Definition} \label {1.1} Let $A=\sum\{A_{ij}\mid i, j\in I\}$ and
$B=\sum \{B_{ij}\mid i, j\in I\}$ are g.m.rings. If $\psi$ is a
map from $A_{ij}$ to $B_{ij}$ for any $i, j\in I$ and satisfies
the following conditions. \begin {eqnarray*} \psi(x+y)&=& \psi
(x)+ \psi (y) \\
\psi (xz)&=& \psi (x) \psi (z) \hbox { for any } x, y\in A_{ij},
z\in A_{j k}, \end {eqnarray*} then $\psi $ is called a
g.m.homomorphism from $A$ to $B$. If $\psi $ is a surjection then
we write $A\cong B$. If $\psi$ is a bijection,  then $\psi$ is
called a g.m.isomorphism, written $A\cong B$.
\end {Definition}

Obviously,  if $\psi$ is a g.m.homomorphism from g.m.ring $A$ to
g.m.ring $B$,  then $\psi$ is a homomorphism from ring $A$ to ring
$B$. Conversely,  it doesn't hold.
\begin {Definition} \label  {1.2}. Let $A=\sum\{A_{ij}\mid i, j\in I\}$ be a g.m.ring
and $B=\sum\{B_{ij}\mid i, j\in I\}$ be a nonempty subset of $A$.
If $B$ is closed under the addition and multiplication,  then $B$
is called a g.m.subring of $A$.

If $B_{ij} A_{j k}\subseteq B_{i k},  A_{ij} B_{j k}\subseteq B_{i
k}$ and $(B_{ij}, +)$ is a subgroup of $(A_{ij}, +)$ for any $i,
j, k\in I$,  then $B$ is called a g.m.ideal of $A$.
\end {Definition}

Obviously,  if $B$ is a g.m.ideal of a g.m.ring $A$,  then $B$ is
an ideal of ring $A$. Conversely,  it doesn't hold.
\begin {Proposition} \label {1.1} Let $B=\sum\{B_{ij}\mid i, j\in I\}$ be a g.m.ideal
of g.m.ring $A$. If let $$A//B=\sum\{ A_{ij}/ B_{ij}\mid i, j\in
I\}$$and define
$$(x+ B_{ij})+(y+ B_{ij})=(x+y)+ B_{ij}$$

$$(x+ B_{ij})(z+ B_{jk})=xz+ B_{ik}$$

for any $x, y\in A_{j k}$,  $z\in A_{jk}$,  $i, j, k\in I$ then $
A//B $ is a g.m.ring.
\end {Proposition}
\begin {Proposition} \label {1.2} If $A\stackrel{\psi}{\cong} B$,  then $A//ker(\psi)\cong B$.
\end {Proposition}
\begin {Proposition} \label  {1.3} If $B$ and $C$ are g.m.ideals of $A$,  then
$(B+C)//C\cong B//(B\cap C)$.
\end {Proposition}
\begin {Proposition} \label  {1.4} If $C$ and $B$ are g.m.ideals of $A$ and  $C\subseteq B$,
then $ A//B\cong(A//C)//(B//C)$.
\end {Proposition}
\begin {Proposition} \label  {1.5}Let $s, t\in I$ and $B_{st}$ be a nonempty subset of
$A_{st}$. If let $D_{ij}= A_{is} B_{st} A_{tj}$ for any $i, j\in
I$,  then $D=\sum\{D_{ij}\mid i, j\in I\}$ is a g.m.ideal of $A$.
\end {Proposition}
\begin {Proposition}  \label {1.6}
If $A^{(t)}$ is a g.m.subring of $A$ for any $t\in \Omega$,  then
$$ (i) \sum _t\{A^{(t)}\mid t\in\Omega \}=\sum \{\sum _tA_{ij}^{(t)}\mid
t\in \Omega \}\mid i, j\in I\}$$

$$(ii)\cap\{A^{(t)}\mid t\in \Omega \}=\sum\{ \cap A_{ij}\mid
t\in \Omega \}\mid i, j\in I\}$$
\end {Proposition}
\begin {Proposition}  \label {1.7}Let $A\stackrel{\psi}{\cong}B$, ${\cal E}_1=\{C\mid C
\hbox { is a g.m.ideal of } A \hbox { and } C \supseteq ker\psi\}$
and ${\cal E } _2=\{D\mid D \hbox { is a g.m.ideal of } B\}$. If
we define $\mu(C)=\psi(C)$ for any $C\in {\cal E}_1$,  then $\mu$
is a bijection from ${\cal E}_1$ into ${\cal E}_2$.
\end {Proposition}
\section {The General Theory of Radicals for G.M.Ring.}

We easily know that the class of all g.m.rings is closed under
g.m.homomorphisms and g.m.ideals. Then we can build up the general
theory of radicals in the class. In this section,  we omit all the
proofs because they are similar to the proofs in ring theory.

\begin {Definition} \label  {2.1} Let $r$ be a property of g.m.rings.
If g.m.ring $A$ is of property $r$ then $A$ is called a
r-g.m.ring. A g.m.ideal $B$ of $A$ is said to be a r-g.m.ideal if
the g.m.ring $B$ is a r-g.m.ring.
\end {Definition}

The property $r$ is called a radical property of g.m.rings if $r$
satisfies the following conditions

(R1)Every g.m.homomorphic image $A'$ of $r$-g.m.ring $A$ is again
a $r$-g.m.ring.

(R2)Every g.m.ring $A$ has a $r$-g.m.ideal,  which contains every
other $r$-g.m.ideal of $A$.

(R3)$A//N$ doesn't contain any non-zero $r$-g.m.ideal. We call the
maximal $r$-g.m.ideal of $A$ the radical of $A$,  written $r(A)$.

If $r(A)=0$, $A$ is called $r$-semisimple.
\begin {Proposition} \label  {2.1}If $r$ is a radical property of g.m.rings,  then
$A$ is a $r$-g.m.ring if and only if $A$ can not be g.m.homomorphically mapped onto
a non-zero $r$-semisimple g.m.ring.
\end {Proposition}
\begin {Proposition} \label  {2.2}Let $r$ is a radical property of g.m.rings. If $B$
is a g.m.ideal of $A$ and $r(A//B)=0$,  then $B\supseteq r(A)$.
\end {Proposition}
\begin {Theorem}  \label {2.1}$r$ is a radical property of g.m.rings if and only if
the following conditions are satisfied

(R1$'$) Every g.m.homomorphic image $A'$ of $r$-g.m.ring $A$ is a
$r$-g.m.ring.

(R2$'$) If every non-zero g.m.homomorphic image of g.m.ring $A$
contains a non-zero $r$-g.m.ideal then $A$ is also an
$r$-g.m.ring.
\end {Theorem}
\begin {Proposition} \label  {2.3}If $r$ is a radical property of g.m.rings,  then
${\cal K}=\{A\mid A \hbox { is an } r \hbox {-semisimple }$
g.m.ring  $ \}$ satisfies the following conditions.

(Q1)Every non-zero g.m.ideal of any g.m.ring from ${\cal K}$ can
be g.m.homomorphically mapped onto a non-zero g.m.ring from ${\cal
K}$.

(Q2)If every non-zero g.m.ideal of g.m.ring $A$ can be
g.m.homomorphically mapped onto a non-zero g.m.ring from ${\cal
K}$,  then $A\in {\cal K}$.
\end {Proposition}
\begin {Proposition} \label  {2.4}If ${\cal K}$ is a class of g.m.rings satisfying conditions
(Q1) and (Q2),  then there exists a radical property $r$ of
g.m.rings such that
$${\cal K}=\{A\mid A \hbox { is a } r\hbox {-semisimple g.m.ring } \}.$$
\end {Proposition}
\begin {Theorem} \label  {2.2}If ${\cal K}$ is a class of g.m.rings satisfying conditions
(Q1) and $\bar{{\cal K}}=\{A\mid \hbox {every non-zero g.m.ideal
of } A \hbox { can be g.m.homomorphically mapped onto a non-zero
}$ g.m.ring in ${\cal K}\}$,  then $\bar{{\cal K}}$ satisfies the
conditions (Q1) and (Q2).
\end {Theorem}
\begin {Definition} \label  {2.2}If ${\cal K}$ satisfies the condition (Q1),  then
${\cal K}$ can determine a radical property $r$ of g.m.rings by
Proposition 2.4 and Theorem 2.2. We call the radical property the
upper radical determined by ${\cal K}\}$,  written $r^{\cal K}$.
\end {Definition}
\begin {Definition} \label  {2.3}Let $r$ be a radical property of g.m.rings and
$$S(r)=\{A\mid A \hbox { is an } r\hbox {-semisimple g.m.ring }\}$$
$$R(r)=\{A\mid A  \hbox { is an } r \hbox {-g.m.ring }\}.$$ If $r'$ is also a radical
property of g.m.rings and $R(r)\subseteq R(r')$ then we say $r\le
r'$.
\end {Definition}
\begin {Theorem}\label   {2.3}Let ${\cal K}$ be a class of g.m.rings satisfying
the condition (Q1) and $r$ be a radical property of g.m.rings. If
${\cal K}\subseteq S(r)$ then $r\le r^{\cal K}$
\end {Theorem}
\begin {Theorem} \label  {2.4}Let $r$ is a radical property of rings. If we define a
property g.m.r of g.m.rings as follows.

A g.m.ring $A$ is a g.m.$r$-g.m.ring if and only if $A$ is a
$r$-ring. Then g.m.r is a radical property of g.m.rings.
Furthermore,  g.m.$r(A)$ is the maximal $r$-g.m.ideal of $A$ for
any g.m.ring $A$.
\end {Theorem}

If we let $r_b, r_k, r_1, r_n, r_j, r_{bm}$, respectively,  denote
the Baer radical,  the nil radical,  the locally nilpotent
radical,  the Neumann regular radical,  the Jacobson radical,
 the Brown McCoy radical of rings,  then following Theorem 2.4,  we can obtain the below
 radical properties of g.m.rings.

g.m.$r_b$, g.m.$r_k$,  g.m.$r_1$,  g.m.$r_n$,  g.m.$r_j$,
g.m.$r_{b m}$.
\section {The Special Radicals of G.M.Ring.}

In this section,  $A$ denotes a g.m.ring. $r_b(A_{ij})$ denotes
the Baer radical of $ A_{ji}$-ring $A_{ij}$.

Shaoxue Liu obtained that if $A$ is an additive category,  the $(
r_b(A))_{ii}= r_b(A_{ii})$ for any $i \in I$ in [3].

In this section,  we obtain that if ${\cal K}$ is a g.m.weakly
special class,  then
$$r^{\cal K}(A)=\cap\{B\mid B \hbox { is a g.m.ideal of } A \hbox { and } A//B \in {\cal K}\}.$$
 Meantime,  we obtain the below conclusion $$r_b(A)=g.m.r_b(A)= \bar r_b(A)=\sum\{ r_b(A_{ij})\mid i, j\in I\}.$$
\begin {Definition} \label  {3.1}A class ${\cal K}$ of g.m.rings is called a g.m.weakly
special class if the following conditions are satisfied

(WS1)Every g.m.ring from ${\cal K}$ is semiprime.

(WS2)Every g.m.ideal of g.m.ring from ${\cal K}$ belongs to ${\cal
K}$.

(WS3)If a g.m.ideal $B$ of g.m.ring $A$ is contained in ${\cal
K}$,  then $A//B^*\in {\cal K}$,  where $B^*=\{x\in A\mid
xA=Ax=0\}$.
\end {Definition}
\begin {Definition} \label  {3.2} A class ${\cal K}$ of g.m.rings is called a g.m.special class
if the following conditions are satisfied

(S1) Every g.m.ring from ${\cal K}$ is prime.

(S2) Every g.m.ideal of g.m.ring from ${\cal K}$ belongs to ${\cal
K}$.

(S3)If a g.m.ideal $B$of g.m.ring $A$ is contained in ${\cal K}$,
then $A//B^*\in {\cal K}$,  where $B^*=\{x\in A\mid xA=Ax=0\}$.
\end {Definition}
\begin {Theorem}  \label {3.1}If ${\cal K}$ is a g.m.weakly special class then
$ r^{\cal K}=\cap\{B\mid B \hbox { is a g.m.ideal of } A $\ \ \ \
$\hbox { and } A//B \in {\cal K}$
\end {Theorem}
{\bf Proof.}Let $r= r^{\cal K}$ and $T=\cap \{B\mid B \hbox { is a
g.m.ideal of } A \hbox { and } A//B \in {\cal K}\}$.We easily show
$T\supseteq r(A)$ by Proposition 2.2.Now we only need show $T$ is
a $r$-g.m.ideal. If $T$ isn't any $r$-g.m.ideal,  then there
exists a g.m.ideal $D$ of $A$ such that $0\ne A//D\in {\cal
K}$.This contradicts the definition of $T$. $\Box$
\begin {Theorem} \label  {3.2} If ${\cal K}$ is a (weakly) special class of rings and let
$g.m.{\cal K}=\{A\mid A  \hbox { is g.m.ring}$ \ \ \ \    and
$A\in {\cal K}$, then g.m.${\cal K}$ is a g.m.(weakly) special
class.
\end {Theorem}

If let $r$ denote $r^{\cal K}$,  then $\bar r$ denotes the upper
radical determined by g.m. ${\cal K}$ of g.m.rings. Following
Theorem 3.2, we can obtain the below radical properties of
g.m.rings.
$$\bar r_b,  \bar r_k,  \bar r_1,  \bar r_j,  \bar r_{bm},  \bar r_n. $$
Following Theorem 2.4 and Theorem 3.2,  we can directly show the
below theorem.
\begin {Theorem}\label   {3.3} Let $r$ be a weakly special radical of rings. Then $$g.m.
r(A)\subseteq r(A)\subseteq \bar r(A) \hbox {  for any g.m.ring } A.$$
$m$-nilpotent element and $m$-sequence have been defined in [5].
Now we use these concepts.
\end {Theorem}
\begin {Theorem} \label  {3.4} $\bar r_b(A)=g.m.r_b(A)=r_b(A)$
\end {Theorem}
{\bf Proof.} Considering Theorem 3.3 we only need show $\bar
r_b(A)$ is a g.m.m-nilpotent ideal of $A$. If $\bar r_b(A)$ isn't
$m$-nilpotent ideal of $A$,  then there exists $x\in \bar r_b(A)$
such that $x$ isn't $m$-nilpotent. That is,
 there exists a
$m$-sequence $\{a_n\}$ in $A$ such that $a_1=x, a_2=a_1u_1a_1,
\cdots ,   $ and $a_n\not=0$ for $n=1, 2, \cdots ,  $ Let ${\cal
E}=\{C\mid C \hbox { is a g.m.ideal of } A \hbox { and }
C\supseteq g.m.r_b(A), C\cap \{a_n\}=\emptyset \}$. By Zorn's
Lemma, there exists a maximal element $V$ in ${\cal E}$. We easily
show that $V$ is a prime g.m.ideal of $A$. If $G//V$ is a non-zero
g.m.$m$-nilpotent ideal of $A//V$,  then there exists $n\in {\bf
Z}^+$ such that $a_n\in G$. Since $G//V$ is $m$-nilpotent, there
exists $m\in {\bf Z}^+$ such that $a_{m+n}\in V$. We get a
contradiction. Then $A//V$ hasn't any non-zero g.m.$m$-nilpotent
ideal. This contradicts $x\in \bar r_b(A)$. $\Box$
\begin {Proposition} \label  {3.1}If $A$ is a semiprime g.m.ring,  then $A_{st}$ is a
semiprime $A_{ts}$-ring for any $s, t\in I$.
\end {Proposition}
{\bf Proof.} If $A_{st}$ is not any semiprime $A_{ts}$-ring, then
there exists a non-zero ideal $B_{st}$ of $A_{st}$ such that
$B_{st}A_{ts}B_{st}=0$. Let $D_{ij}=A_{is}B_{st}A_{tj}$ for any
$i, j\in I$. Since $D$ is an ideal of $Aa$ by Proposition 1.5 and
$D_{ik}D_{kj}=A_{is}B_{st}A_{tk}A_{ks}B_{st}A_{tj} \subseteq
A_{is}B_{st}A_{ts}B_{st}A_{tj}=0$ for any $i, j, k\in I,
DD=0.$Since $A$ is a semiprime ring, $D=0$, i.e.$A_{i s}B_{st}A_{t
j}=0$ for any $i, j\in I$ and $A(B_{st}E(s, t))A=0$. Then
$B_{st}=0$. We get a contradiction. $\Box$
\begin {Proposition} \label  {3.2}If $A$ is a prime g.m.ring,  then $A_ {ij}$ is a prime
$A_ {ji}$-ring for any $i, j\in I$
\end {Proposition}
\begin {Theorem} \label  {3.5}$\bar r_b(A)\supseteq\sum\{r_b(A_ {ij})\mid i, j\in I\}.$
\end {Theorem}
{\bf Proof.} $\bar r_b(A)=\cap\{B\mid B \hbox { is a prime
g.m.ideal of } A\}=\sum\{\cap\{ B_{ij}\mid B \hbox { is a prime
g.m.ideal}$ of $ A\} \mid i, j\in I\}$ by Proposition 1.6. If
$A//B$ is a prime g.m.ring,  then $A_{ij}//B_{ij}$ is  a prime $
A_{ji}/B_{ji}$-ring by Proposition 3.2, we easily show that
$A_{ij}/B_{ij}$ is also a prime  $A_{ji}$-ring. Then
$r_b(A_{ij}/B_{ij})=0$. Considering [5 Theorem 3.1] and [5
Proposition 1.5],  we have
$$r_b(A_{ij})\subseteq B_{ij} \hbox { for any } i, j\in I.$$ Therefore $\bar
r_b(A)\supseteq\sum\{r_b(A_{ij})\mid i, j\in I\}$. $\Box$
\begin {Theorem} \label  {3.6}$g.m.r_b(A)\subseteq\sum\{ r_b(A_{ij}) \mid i, j\in
I\}$
\end {Theorem}
{\bf Proof.}We first show that if $r_b(A)=A$ then
$r_b(A)=\sum\{r_b(A_{ij})\mid i, j\in I\}$. Let $x\in A_{ij}$. If
$\{ a_s\mid s=1, 2\cdots   \}$ is an $m$-sequence with $a_1=x$,
i.e there exists $u_s\in A_{ji}$ such that $a_{s+1}=a_su_sa_s$ for
$ s=1, 2\cdots   $. Let $b_s=a_sE{{i, j}}$, $v_s=u_sE{{i, j}}$ for
$ s=1, 2\cdots  $. Since $\{b_s\mid s=1, 2\cdots  \}$ is an
$m$-sequence in ring $A$ and $r_b(A)=A$,  there exists $k\in {\bf
Z}^+$ such that $b_k=0$,  i.e. $a_k=0$ and $x$ is $m$-nilpotent.
Therefore $r_b(A_{ij})=A_{ij}$.

Next we show $g.m.r_b(A)\subseteq \sum\{ r_b(A_{ij})\mid i, j\in
I\}$. Let $N=g.m. r_b(A)$. Since $ r_b(N)=N$, $N_{ij}$ is a
$r_b$-$N_{ji}$-ring by the above proof. Now we show that $N_{ij}$
is $r_b$-ideal of $A_{ji}$-ring $A_{ij}$. Let $a\in N_{ij}$ and
$\{a_n\}$ be an $m$-sequence in $A_{ji}$-ring $N_{ij}$,  where
$a_1=a, a_{n+1}=a_nu_na_n, u_n\in A_{ji}$ for $n=1, 2, \cdots ,  $

Let $b_n=a_{2n}$. Since $a_{2n+1}=
a_{2n}u_{2n}a_{2n}u_{2n+1}a_{2n}u_{2n}a_{2n}=
a_{2n}v_na_{2n}=b_nv_nb_n$,
 where $v_n=u_{2n}a_{2n}u_{2n+1}a_{2n}u_{2n}\in N_{ji}$, $\{b_n\}$ is an $m$-sequence in $N_{ji}$-ring $N_{ij}$.
 Therefore there exists $m\in {\bf Z}^+$ such that $b_m=0$. Thas is , $\{a_n\}$ is
 an
 $m$-nilpotent sequence. Thus $N_{ji}$ is $r_b$-ideal of $A_{ji}$-ring $A_{ij}$, i.e
 $ r_b(A_{ij})\supseteq N_{ij} \hbox { for any } i, j\in I$, and $g.m.r_b(A)\subseteq \sum\{ r_b(A_{ij})
 \mid i, j\in I\}$.$\Box$

 Following Theorem 3.4, Theorem 3.5 and Theorem 3.6,  we
have the below conclusion.
\begin {Theorem} \label  {3.7}$r_b(A)=\bar r_b(A)=g.m.r_b(A)=\sum\{ r_b(A_{ij})
\mid i, j\in I\}$
\end {Theorem}
\begin {Proposition}\label   {3.3}Let $A$ be a ring and $M$ a $\Gamma$-ring with
$M=\Gamma=A$. Then $r_b(A)=r_b(M)$.
\end {Proposition}
{\bf Proof.}Let $W(A)=\{x\in A\mid x \hbox { is an } m\hbox {
-nilpotent element of ring } A\}$. Let $W(M)=\{x\in M\mid x \hbox
{ is a} m\hbox { -nilpotent element of }\Gamma \hbox {-ring }
M\}$. Considering \cite  [Theorem 3.9 and Definition 3.4] {ZC91},
we have $W(A)=r_b(A)$, $W(M)=r_b(M)$ and $W(A)=W(M)$. Then
$r_b(A)=r_b(M)$.

Following Proposition 3.3,  we have $r_b(A_{ii})$ is also the Baer
radical of ring $A_{i i}$ for any $i\in I$.

\end {document}